\documentclass[12pt]{article}
\usepackage{amsmath,amssymb,amsthm}
\usepackage{mathrsfs}
\usepackage{authblk}

\newcommand{\handout}[5]{
  \noindent
  \begin{center}
  \framebox{
    \vbox{
      \hbox to 5.78in {\hfill #2 }
      \vspace{4mm}
      \hbox to 5.78in { {\Large \hfill #5  \hfill} }
      \vspace{2mm}
      \hbox to 5.78in { {\em #3 \hfill #4} }
    }
  }
  \end{center}
  \vspace*{4mm}
}

\newtheorem{theorem}{Theorem}
\newtheorem{corollary}[theorem]{Corollary}
\newtheorem{lemma}[theorem]{Lemma}

\newtheorem{proposition}[theorem]{Proposition}
\theoremstyle{definition} 
\newtheorem{definition}[theorem]{Definition}

\newtheorem{remark}[theorem]{Remark}
\newtheorem{example}[theorem]{Example}

\advance \topmargin by -\headheight
\advance \topmargin by -\headsep
\evensidemargin \oddsidemargin
\usepackage{xcolor}
\usepackage[bookmarks=true]{hyperref}
\hypersetup{
	colorlinks,
	linkcolor={red!70!black},
	citecolor={red!90!black},
	urlcolor={green!80!black},
	linktoc=all, 
}

\usepackage[backend=bibtex,style=alphabetic,sorting=nyt,maxbibnames=99,giveninits=true]{biblatex}  
\renewbibmacro{in:}{} 
\DeclareFieldFormat{pages}{#1}
\DeclareFieldFormat[article]{number}{ no. #1}
\bibliography{Bbib} 

\usepackage{booktabs}

\begin{document}

\title {Vandermonde sets and hyperovals}
\author[ ]{Kanat Abdukhalikov and Duy Ho} 
\affil[ ]{UAE University, PO Box 15551, Al Ain, United Arab Emirates}
\affil[ ]{\textit {abdukhalik@uaeu.ac.ae, duyho92@gmail.com}}

\date{ }
\maketitle

\begin{abstract}
	We consider relationships between Vandermonde sets  and hyperovals. 
	Hyperovals are Vandermonde sets, but, in general, Vandermonde sets are not hyperovals. 
	We give necessary and sufficient conditions for a  Vandermonde set to be a hyperoval.    
	Therefore, we provide purely algebraic criteria for existence of hyperovals. 
	Furthermore, we give necessary and sufficient conditions for the existence of hyperovals in terms of $g$-functions, 
	which can be considered as an analog of Glynn's Theorem for o-polynomials.  
\end{abstract}

\section{Introduction}

A finite set of $t$ elements in a finite field is called a Vandermonde set if all power sums of degrees up to $t-2$ 
of these elements are equal to $0$. 
We consider relationships between Vandermonde sets  \cite{gacs2003,sziklai2008} and hyperovals. 
Hyperovals are Vandermonde sets, but, in general, Vandermonde sets are not hyperovals. 
We found an additional set of power sums that their equality to $0$ provides necessary and sufficient 
conditions for a  Vandermonde set to be a hyperoval.    
Therefore, we provide purely algebraic criteria for existence of hyperovals in terms of power sums. 

Usually hyperovals are described in terms of $o$-polynomials, and an alternative description of hyperovals 
in terms of special $g$-functions are given in \cite{kanat2017,kanat2019}. 
We note that $g$-function approach is more efficient in cases when o-polynomial looks very complicated,  
in particular for Subiaco, Adelaide, Lunelli-Sce and sporadic O'Keefe-Penttila hyperovals \cite{kanat2019b,kanat2019c}.  
We give necessary and sufficient conditions for the existence of hyperovals in terms of $g$-functions, 
which can be considered as an analog of Glynn's Theorem \cite{glynn1989} for $o$-polynomials.  
Provided conditions for existence of hyperovals are suitable for efficient computer search. 

The paper is organized as follows.  
We  first recall in Section 2 preliminary results on projective planes, affine planes and hyperovals. 
In Section 3, we show in Theorem \ref{thmhovalvset} that hyperovals are Vandermonde sets.
In Section 4, we provide three new characterizations of a hyperoval via power sums  of its points. These are summarized in Theorem \ref{thmbracket} and  Theorem \ref{thmhovalD2}. We also describe examples of Vandermonde sets that are not hyperovals.
In Section 5, we describe conditions on $g$-functions and $\rho$-polynomials for the existence of hyperovals in  Theorem \ref{thmcoefg} and Theorem \ref{thmcoeff}.
In Section 6, we consider the Gram matrix for the elements of a hyperoval.

\section{Preliminaries}

\subsection{Polar representation}
In our paper we consider finite fields of characteristics $2$ only. 
Let $F = \mathbb{F}_{2^m}$ be a finite field of order $q = 2^m$. Consider $F$ as a subfield of $K = \mathbb{F}_{2^n}$, where
$n = 2m$, so $K$ is a two dimensional vector space over $F$.

The \textit{conjugate} of $x \in K$ over $F$ is $$\bar{x} = x^q.$$
Then the \textit{trace} and the \textit{norm} maps from $K$ to $F$ are
$$T(x) = Tr_{K/F} = x + \bar{x} = x + x^q,$$
$$N(x) = N_{K/F}(x) = x\bar{x} = x^{1+q}.$$

The \textit{unit circle} of $K$ is the set of elements of norm $1$:
$$S = \{u \in K \mid u \bar{u}= 1 \}.$$

Therefore, $S$ is the multiplicative group of $(q +1)$st roots of unity in $K$. 
Since $F \cap S = \{1\}$, each non-zero element of $K$ has a unique polar coordinate representation $x = \lambda u$
with $\lambda \in F^*$ and $u \in S$. For any $x \in K^*$ we have 
$\lambda = \sqrt{x \bar{x}}$ and $u = \sqrt{x /\bar{x}}$.

One can define nondegenerate bilinear form $\langle \cdot, \cdot \rangle : K \times K \rightarrow F$ by
$$\langle x, y \rangle = T(x\bar{y}) = x\bar{y} + \bar{x}y.$$
Then the form  $\langle \cdot, \cdot \rangle$ is alternating and symmetric, that is, 
$\langle a, a \rangle=0$ and 
$\langle a, b \rangle= \langle b, a \rangle$.

%
\subsection{Affine and projective planes}
Consider points of a projective plane $PG(2, q)$ in homogeneous coordinates as triples  $(x:y:z)$, where $x,y,z\in F$, $(x,y, z) \ne (0, 0, 0)$, and we identify $(x : y : z)$ with $(\lambda x : \lambda y : \lambda z)$, $\lambda \in F$. Then points of $PG(2, q)$ are
$$
\{(x : y : 1) \mid x \in F, y \in F \} \cup \{ (x : 1 : 0) \mid x \in F\} \cup \{(1 : 0 : 0)\}.
$$
For $a, b, c \in F$, $(a, b, c) \ne (0, 0, 0)$, the line $[a : b : c]$ in $PG(2, q)$ is defined as
$$
[a : b : c] = \{(x : y : z) \in PG(2, q) \mid ax + by + cz = 0\}.
$$

Triples $[a : b : c]$ and $[\lambda a : \lambda b : \lambda c]$ with $\lambda \in F^*$ define same lines. The point $(x : y : z)$ is incident with the line $[a : b : c]$ if and only if $ax + by + cz = 0$. We shall call points of the form $(x : y : 0)$ the points at infinity. Then $[0 : 0 : 1]$ indicates the line at infinity.

We define an affine plane $AG(2, q) = PG(2, q) \backslash [0 : 0 : 1]$, so points of this affine plane $AG(2, q)$ are $\{(x : y : 1) \mid  x, y \in F \}$.
Associating $(x : y : 1)$ with $(x, y)$ we can identify points of the affine plane $AG(2, q)$ with elements of the vector space $V(2, q) = \{(x, y) \mid x, y \in F\}$, and we will write $AG(2, q) = V (2, q)$.
Lines in $AG(2, q) = V (2, q)$ are $\{(c, y) \mid y \in F\}$ and $\{(x, xb + a) \mid x \in F\}, a, b, c \in F$.
These lines can be described by equations $x = c$ and $y = xb + a$.

We introduce now other representation of $PG(2, q)$ using the field $K$. Consider pairs $(x : z)$, where $x \in K, z \in F$, $x \ne 0$ or $z \ne 0$, and we identify $(x : z)$ with $(\lambda x : \lambda z),  \lambda \in F^*$. 
Then points of $PG(2, q)$ are
$$
\{(x : 1) \mid x \in K\} \cup \{(u : 0) \mid u \in S \}.
$$
For $\alpha\in K$ and $\beta\in F$ we define lines $[\alpha:\beta]$ in $PG(2, q)$ as
$$
[\alpha:\beta] = \{(x:z) \in PG(2,q) \mid \langle\alpha,x\rangle+\beta z = 0\}.
$$
Pairs $[\alpha:\beta]$ and $[\lambda\alpha:\lambda\beta]$ with $\lambda \in F^*$ define same lines. The point $(x : z)$ is incident with the line $[\alpha:\beta]$ if and only if $\langle\alpha,x\rangle+\beta z = 0$. The element $u_\infty = (u : 0), u \in S$, will be referred to as the point at infinity in the direction of $u$. So $[0 : 1]$ indicates the line at infinity.

We define an affine plane $AG(2, q) = PG(2, q) \backslash [0 : 1]$, so points of  this affine plane $AG(2,q)$ are $\{(x : 1) \mid x \in K\}$.
Associating $(x : 1)$ with $x \in K$ we can identify points of the affine plane $AG(2,q)$ with elements of the field $K$, and we write $AG(2, q) = K.$ 
Lines of $AG(2, q)=K$ are of the form
$$
L(u,\mu) = \{x \in K \mid \langle u, x\rangle + \mu = 0\},
$$
where $u \in S$ and $\mu \in F$ (cp. \cite[subsection 2.1]{ball1999}).

Throughout the paper, we will consider these two representations of the projective plane $PG(2, q)$, and for each of such projective planes we consider a fixed affine plane $AG(A, q)$ described above. 
They will be written as $AG(2, q) = V (2, q)$ and $AG(2, q) = K$.

\subsection{Hyperovals and representations}
In the projective plane $PG(2, q)$, $q=2^m$, 
an \textit{oval} is a set of $q +1$ points, no three of which are collinear. Any line of the plane meets the oval $\mathcal{O}$ at either $0, 1$ or $2$ points and is called exterior, tangent or secant, respectively. All the tangent lines to the oval $\mathcal{O}$ concur at the same point $N$, called the \textit{nucleus} of $\mathcal{O}$. The set $\mathcal{H} = \mathcal{O} \cup N$ becomes a \textit{hyperoval}, that is a set of $q + 2$ points, no three of which are collinear. Conversely, by removing any point from hyperoval one gets an oval.

By the fundamental theorem of projective geometry, any hyperoval of $PG(2, q)$  is equivalent to a hyperoval $\mathcal{H}$  containing the points $(1: 0: 0), (0: 0: 1), (0: 1: 0)$ and $(1: 1: 1)$. 
Consequently we may write $\mathcal{H}$ in the form
$$
\mathcal{H}= \{(t:f(t):1) \mid  t \in F \} \cup \{(1: 0:0), (0: 1: 0) \},
$$
where $f$ induces a permutation of $F$ such that $f(0) = 0$ and $f(1) = 1$.
By applying Lagrange's Interpolation Formula, it can be shown (cp. \cite{hirschfeld1998}) that $f$ can be expressed uniquely as a polynomial of degree $q - 2$ over the field $F$. 
The permutation polynomials which arise from hyperovals in this way are called \textit{o-polynomials} (cp. \cite{cherowitzo1986}, \cite{cherowitzo1996}).

%


As described in \cite{kanat2019}, the hyperoval $\mathcal{H}$ can also be represented in $PG(2, q)$ using
the field $K$ as
$$
\mathcal{H}=\left\{ \dfrac{u}{g(u)} \mid u  \in S\right\} \cup \{0\}, 
$$
for some function $g: S \rightarrow F$.
If $g(u) = 0$ then we assume that $u/g(u) = u_{\infty}$ is the element
at infinity in the direction $u$. 
Such a function $g : S \rightarrow F$ is said to be a \textit{$g$-function}. 
We can assume $g(u)\ne 0$ for all $u \in S$, by taking in place of $g(u)$ an equivalent function $g(u)+\langle c,u \rangle$ with appropriate $c \in K$, cp. \cite{kanat2017}. Furthermore, we can assume that
$$
g(u)=a_0+a_1u+ \cdots a_qu^q,
$$
where $a_0 \in \mathbb{F}_2$, $a_i \in K$  and $a_{q+1-i}=a_i^q$  for $1 \le i \le q/2$.
%
%

In thesis \cite{deorsey2015}, following ideas from \cite{fisher2006}, the \textit{$\rho$-polynomials} were introduced.
We note that the $\rho$-polynomials and $g$-functions are connected in the following way: $g(u) = 1/\rho(u)$.

\section{Hyperovals are Vandermonde sets}
\begin{definition} \theoremstyle{definition} Let $1 < t < q$. A set  $T := \{y_1, \cdots, y_t\} \subseteq GF(q)$ is a \textit{Vandermonde set} if
$$
\pi_k(T):=\sum_{i=1}^{t} y_i^k=0,
$$
for all $1 \le k \le t-2$ (cp. \cite{gacs2003, sziklai2008}). The set $T$ is a \textit{super-Vandermonde set} if it is a Vandermonde set and $\pi_{t-1}(T)=0$.
\end{definition}

\begin{lemma}[cp. \cite{sziklai2008}] \label{Vandermondetranslate} The Vandermonde property is invariant under transformations $y \rightarrow ay+b$, $(a \ne 0)$ if and only if $t$ is even or $b=0$.

\begin{proof} 
	We have	$\pi_k(T):= \sum_i y_i^k =0$ for all $1 \le k \le t-2$.  Denote $T'$ the transformed $T$. For $1 \le k \le t-2$,
	\begin{align*}
	\pi_k(T') &=\sum_{i=1}^t (ay_i+b)^k  \\
	&= a^k \left(\sum_i y_i^k \right) +tb^k + \sum_{i=1}^t\sum_{j=1}^{k-1}  \binom{k}{j}a^jy_i^jb^{k-j} \\
	&=tb^k + \sum_{j=1}^{k-1} \left( \binom{k}{j}a^jb^{k-j}\sum_{i=1}^t y_i^j \right) \\
	&= tb^k.
	\end{align*}
	The proof now follows. 
\end{proof}
\end{lemma}

\begin{corollary}[cp. \cite{sziklai2008}] If $T$ is a Vandermonde set containing the zero element, then $T \backslash \{ 0 \}$ is a super-Vandermonde set.   In particular, if $T$ is a Vandermonde set and $t$ is even, then for any $a \in T$, the translate $T-a$ is a Vandermonde set containing the zero element. 
\end{corollary}

\begin{lemma} \label{Voval} Let $O=\{y_1, \cdots y_{q+1} \}$ be an oval with points in $GF(q^2)$ and nucleus at $0$. Then $O$ is a super-Vandermonde set of $q+1$ points. 

\begin{proof} The set $H:= O \cup \{ 0 \}$ is a hyperoval. Let $x \not \in  H$. Since $(x-y_i)^{q-1}$  represents the slope of line going through $x$ and $y_i$ (cp.  \cite{ball1999}) and since every line going through $x$ intersects $H$ at either $0$ or $2$ points, we have
	\begin{align*} 
	x^{q-1}+\sum_{i=1}^{q+1}(x-y_i)^{q-1} & =0. 
	\end{align*}
	This implies the polynomial
	\begin{align*} 
	\chi(X) &= X^{q-1}+\sum_{i=1}^{q+1}(X-y_i)^{q-1}
	\end{align*}
	has at least $q^2-q-2$ roots. On the other hand, the degree of $\chi(X)$ is at most $q-2$, so $\chi(X)$ must be the zero polynomial. This implies
	$$
	\pi_k(O) = \sum_{i=1}^{q+1} y_i^k=0, 
	$$
	for $1 \le k \le q-1$. 
	Also,
	
	$$\pi_q(O)= \sum_{i=1}^{q+1} y_i^q = \left( \sum_{i=1}^{q+1} y_i \right)^q = 0.$$
	
	Therefore $O$ is a super-Vandermonde set. \qedhere

\end{proof}
\end{lemma}

\begin{theorem} \label{thmhovalvset} A hyperoval with points in $GF(q^2)$ is a Vandermonde set.  
\begin{proof} Let $N$ be a  hyperoval with points in $GF(q^2)$.  Let $N_0$ be a translation of $N$ containing $0$. Then $O:= N_0 \backslash \{0\}$ is an oval with points in $GF(q^2)$ and nucleus at $0$. By Lemma \ref{Voval}, $O$ is a super-Vandermonde set. Hence $N_0$ is a Vandermonde set. By Lemma \ref{Vandermondetranslate}, $N$ is also a Vandermonde set. 
\end{proof}
\end{theorem}

\section{Hyperovals and power sums of points} 
In this section, we consider the set $H:= \{ u/g(u) \mid u \in S \} \cup \{ 0\}$.  As noted in Section 2, we can assume $g(u) \ne 0$ for all $u \in S$. 
We rewrite elements of $S$ as $S:= \{ u_0, u_1, \cdots u_q \}$. For $1 \le j \le q+1$, let $y_j:=u_j/g(u_j)$. For $1 \le k \le q$, let $\pi_k:= \sum_{j=1}^{q+1} y_j^k$.

\begin{lemma} \label{hovaleven} A set $N$ of $q+2$ points in $GF(q^2)$ is a hyperoval if and only if every line intersects $N$ at an even number of points. 

\begin{proof} If $N$ is a hyperoval then every line intersects $N$ at $0$ or $2$ points. Assume every line intersects $N$ at an even number of points. Let $p$ be a point on $N$. Then the $q+1$ lines going through $p$ each contains at least one more point of $N$. Since there are only $q+1$ points in $H \backslash \{p\}$, each of these lines contain exactly two points of $N$ (including $p$). In particular, no three points of $N$ are collinear and so $N$ is a hyperoval. 
\end{proof}

\end{lemma}

\begin{example} Vandermonde sets of size $q+2$ are not necessarily hyperovals. For $q=8$, let $\lambda,\mu \in F^*$ such that $1+\lambda^3+\mu^3=0$. Let 
$$
H:= \{0, 1, \omega, \bar{\omega}, \lambda, \lambda\omega, \lambda\bar{\omega},\mu, \mu\omega, \mu\bar{\omega} \},
$$
where $\omega \in S$ such that $\omega^3=1$. Then $H$ is a Vandermonde set but not a hyperoval, as the line $\langle 1,x \rangle =0$ intersects $H$ at four points.
\end{example}

So from now on we will concentrate on sets of the form $H:= \{ u/g(u) \mid u \in S \} \cup \{ 0\}$.

\begin{theorem}  \label{thmbracket} The set $H:= \{ u/g(u) \mid u \in S \} \cup \{ 0\}$ is a hyperoval if and only if one of the followings holds. 
\begin{enumerate}
	\item The equation
	\begin{equation} \label{eqnkanat}
	g(u) + \langle u,b\rangle = 0
	\end{equation}
	has an even number of solutions $u \in S$ for each $b \in K$. 
	\item 	For each $v \in S$ and $1 \le k \le q$,
	\begin{equation}  \label{eqglynn}
	\sum_{u \in S} \left\langle v, \dfrac{u}{g(u)}\right\rangle^k=0.
	\end{equation}
	
\end{enumerate}
\begin{proof}
	\begin{enumerate}
		\item  If $H$ is a hyperoval, then by \cite[Theorem 4.4]{kanat2019} the equation \eqref{eqnkanat} has $0$ or $2$ solutions for each $b \in K$. Assume the converse. For $v\in S, \mu \in F$, we want to show the line $L$ defined by
		$$
		\langle v, x \rangle +\mu=0
		$$
		intersects $H$ at an even number of points. There are two cases depending on $\mu.$
		
		\begin{enumerate}
			\item $\mu = 0$. Then $L$ intersects $H$ at two points $0$ and $v/g(v)$. 
			\item $\mu \in F \backslash \{0\}$. Let  $b=\dfrac{v}{\mu}$. We have
			\begin{align*} 
			g(u) = \langle u,b\rangle &\iff g(u)= \langle u,v/ \mu \rangle \\
			&\iff g(u) = 1/\mu \cdot\langle u,v \rangle \\
			&\iff \mu = 1/g(u) \cdot \langle u,v \rangle \\
			&\iff \mu = \langle v,u/g(u) \rangle \\
			&\iff  \langle v,u/g(u) \rangle+\mu=0. 
			\end{align*}
			This implies the line $L$ intersects $H$ at an even number of points. The claim now follows from Lemma \ref{hovaleven}.
		\end{enumerate}  
		\item Assume $H$ is a hyperoval. Fix $v \in S$. For each $\mu \in F$, the line 
		$
		\left\langle v, x \right\rangle = \mu
		$
		intersects $H$ at either $0$ or $2$ points. In particular, the equation
		$$
		\left\langle v, \dfrac{u}{g(u)}\right\rangle = \mu
		$$
		has either $0$ or $2$ solutions $u \in S$. This implies condition \eqref{eqglynn} is true for each $v \in S$ and $1 \le k \le q$.
		
		Conversely, assume condition \eqref{eqglynn} is true for each $v \in S$ and $1 \le k \le q$.  We will show that every line $\langle v, x\rangle = \mu$ intersects $H$ at an even number of points. Fix $v \in S$.   For $\mu=0$, the line  $\langle v, x\rangle =0$ intersects $H$ at $0$ and $v/g(v)$.   For each $\mu \in F \backslash \{0\}$, let 
		$$U_\mu := \left\{ u \in S \backslash \{v\} \mid \left\langle v, \dfrac{u}{g(u)}\right\rangle = \mu \right\}.$$ We note that, for each $1 \le k \le q$,  
		\[
		\sum_{u \in U_\mu} \left\langle v, \dfrac{u}{g(u)}\right\rangle^k= 
		\begin{cases}
		\mu^k & \text{if $|U_\mu|$ is odd, }\\
		0         & \text{if $|U_\mu|$ is even.}
		\end{cases} 
		\]

		Let $\Omega:= \{ \mu \in F \backslash \{0\} \mid   |U_\mu| \text{ is odd} \}$.  Assume that $\Omega \ne \varnothing$. Then
		$$
		\sum_{u \in S} \left\langle v, \dfrac{u}{g(u)} \right\rangle^k
		=\sum_{\mu \in \Omega}	\sum_{u \in U_\mu} \left\langle v, \dfrac{u}{g(u)}\right\rangle^k 
		= \sum_{\mu \in \Omega} \mu^k = 0.
		$$
		Also, the sets $U_\mu$ partition the set $ S \backslash \{v\}$ of size $q$,  so  that $|\Omega|$ is even. In particular, the sum of the vectors $(1, \mu, \mu^2, \cdots, \mu^{q-2})$, $\mu \in \Omega$, is the zero vector. 
		On the other hand, the Vandermonde's determinant implies that these vectors are linearly independent. Hence $\Omega = \varnothing$.

		We have shown that  every line $\langle v, x\rangle = \mu$ intersects $H$ at an even number of points. By Lemma \ref{hovaleven}, $H$ is a hyperoval. \qedhere
	\end{enumerate}
\end{proof}
\end{theorem}

%
%
%
%
%
%
%
\begin{lemma} \label{bracketpower} The power of the bilinear form $\langle \cdot, \cdot \rangle$ is given by 
$$
\langle a,b \rangle^k = \sum_{i=0}^{\lfloor(k-1)/2\rfloor} \binom{k}{i} \langle a^{iq+k-i},b^{iq+k-i} \rangle.
$$

\begin{proof}  
	Assume $k$ is odd. We have
	\begin{align*}
	\langle a,b \rangle^k &= (a^qb+ab^q)^k= \sum_{i=0}^k \binom{k}{i} a^{iq}b^{i}a^{k-i}b^{(k-i)q}= \sum_{i=0}^k \binom{k}{i} a^{iq+k-i}b^{i+(k-i)q} \\
	&= \sum_{i=0}^{(k-1)/2} \binom{k}{i} a^{iq+k-i}b^{i+(k-i)q}
	+ \sum_{i={(k+1)/2}}^{k} \binom{k}{i}  a^{iq+k-i}b^{i+(k-i)q} \\
	&= \sum_{i=0}^{(k-1)/2} \binom{k}{i} a^{iq+k-i}b^{i+(k-i)q}
	+ \sum_{i=0}^{(k-1)/2} \binom{k}{i}  a^{i+(k-i)q}b^{iq+k-i} \\
	&= \sum_{i=0}^{(k-1)/2}\binom{k}{i} \left(  a^{iq+k-i}b^{i+(k-i)q}+ a^{i+(k-i)q}b^{iq+k-i} \right) \\
	&=\sum_{i=0}^{(k-1)/2} \binom{k}{i}
	\langle a^{iq+k-i},b^{iq+k-i} \rangle.
	\end{align*}   
	Similar to the above, when $k$ is even, we have
	\begin{align*}
	\langle a,b \rangle^k  &= \sum_{i=0}^{k/2-1} \binom{k}{i} \langle a^{iq+k-i},b^{iq+k-i} \rangle + \binom{k}{k/2} \langle a^{iq+k-i},b^{iq+k-i} \rangle \\
	&= \sum_{i=0}^{k/2-1} \binom{k}{i}\langle a^{iq+k-i},b^{iq+k-i} \rangle,
	\end{align*}   
	since $\binom{k}{k/2}$ is even.
	\qedhere

\end{proof}
\end{lemma}
To state the next theorem we define the following partial ordering $\preceq$ on the set of nonnegative integers. If
$$
b=\sum_{i=0}^{m}b_i2^i \text{ and } c=\sum_{i=0}^{m}c_i2^i
$$
(where each $b_i$ and each $c_i$ is either $0$ or $1$), then $b \preceq c$ if and only if $b_i \le c_i$ for all $i$ (cp. \cite{glynn1989}, \cite{okeefe1991}). 
In other words, $b \preceq c$ if and only if all nonzero terms appearing in the
binary expansion of $b$ also appear in the binary expansion of $c$.

Let
\begin{equation} \label{eqnM}
\mathscr{M}:=\{ (i,k) \mid  1 \le k \le q-2, 0\le i \le \lfloor(k-1)/2\rfloor, 
i \preceq k \},
\end{equation}
and 
\begin{equation} \label{eqnD}
D:= \{iq+k-i\mid (i,k) \in \mathscr{M} \}. 
\end{equation}
\begin{theorem} \label{hovalD}
The set $H:= \{ u/g(u) \mid u \in S \} \cup \{ 0\}$ is a hyperoval  if and only if $\pi_{d} = 0$ for all $d \in D$. 
\begin{proof} 
	\begin{enumerate}
		\item  We first prove that $H$ is a hyperoval  if and only if $\pi_{d} = 0$ for all $d \in D'$, where $$D':= \{iq+k-i\mid  1\le k\le q, 0\le i \le \lfloor(k-1)/2\rfloor, i \preceq k \}. $$ 
		
		Fix $1\le k\le q$. 	For each $v \in S$,  by Lemma \ref{bracketpower}, we have
		\begin{align*}
		\sum_{u \in S} \left(\langle v, \dfrac{u}{g(u)}\rangle \right)^k 
		&= \sum_{j=1}^{q+1} \sum_{i=0}^{\lfloor(k-1)/2\rfloor} \left( \binom{k}{i} \left\langle v^{iq+k-i},y_j^{iq+k-i} \right\rangle\right)\\
		&= \sum_{i=0}^{\lfloor(k-1)/2\rfloor}\sum_{j=1}^{q+1}\left(\binom{k}{i} \left\langle v^{iq+k-i},y_j^{iq+k-i}\right\rangle\right)\\ 
		&= \sum_{i=0}^{\lfloor(k-1)/2\rfloor} \left(\binom{k}{i} \left\langle v^{iq+k-i},\pi_{iq+k-i}\right\rangle\right)\\
		&= \sum_{\substack{ i \preceq k \\0 \le i \le \lfloor(k-1)/2\rfloor }} \left\langle v^{iq+k-i},\pi_{iq+k-i}\right\rangle = P_k(v),
		\end{align*}
		
		where $P$ is the polynomial obtained from expanding the terms in the last sum. We can assume $P_k$ has degree at most $q$, as $v^{q+1}=1$.
		
		The power of $v$ in $P_k(v)$ is either of the form $k-2i$ or $q+1+2i-k$. If $i \ne i'$, then $k-2i \ne k-2i'$, and  $k-2i \ne 2i'-k+q+1$.	Hence $P_k(v)$ has the form $P_k(v) = \sum \pi_rv^s$. 
		
		By Theorem \ref{thmbracket}, $H$ is a hyperoval if and only if	$P_k(v)=0$ has $q+1$ roots $v \in S$ for each $1 \le k \le q$. Equivalently, the coefficients of $P_k$ are zeros, that is, $\pi_d=0$ for all $d \in D'$.

		\item Since $D \subseteq D'$, if $\pi_{d'}=0$ for all $d' \in D'$, then $\pi_{d}=0$ for all $d \in D$.  Conversely, we let $d'=iq+k-i \in D'$ and consider the two cases $k=q$ and  $k=q-1$. 
		In the case $k=q$,  there is only one possibility $i=0$ and $d'=q$. In particular,  $\pi_{d'}=0$ if and only if $\pi_1=0$.

		The case  $k=q-1$ occurs if and only if  $q-1 \mid d'$. This implies $d'=(q-1)t$ for some $1 \le t \le q$ and
		$$
		\pi_{d'} = \sum_{j=1}^{q+1} (y_j^{q-1})^t. 
		$$
		The elements $y_j^{q-1}$ are in $S$ and pairwise distinct, so that the set $\{ y_j^{q-1} \mid 0 \le j \le q \}$ is $S$ and hence a Vandermonde set. It follows that $\pi_{d'}=0$.

		Therefore, $\pi_{d'}=0$ for all $d' \in D'$ if and only if $\pi_{d}=0$ for all $d \in D$. 
	\end{enumerate}

	The theorem now follows from part 1) and 2). \qedhere

\end{proof}
\end{theorem} 
%
%
%
%
%
%
%
%
%
For $d\in D$ and $l \in \mathbb{Z}$, if $d'\equiv2^ld \pmod{q^2-1}$, then $\pi_d=\pi_{d'}$. And so Theorem \ref{hovalD} can be improved as follows.
Let $\sim$ be the equivalence relation on $D$ defined by $x \sim y$ if and only if there  exists $l \in \mathbb{Z}$ such that $x\equiv 2^ly \pmod{q^2-1}$. Let 
$$\mathcal{D}:=D/ \sim. $$
\begin{theorem} \label{thmhovalD2} 
The set $H:= \{ u/g(u) \mid u \in S \} \cup \{ 0\}$ is a hyperoval  if and only if $\pi_{d} = 0$ for all $d \in \mathcal{D}$. 
\end{theorem}

\begin{remark}For $q =2^m \le 128$,  the set $\mathcal{D}$ is given in Table \ref{tabledset}.
\end{remark}

\begin{remark} For $q=4$  and $q=8$, from Table \ref{tabledset} it follows that $H:= \{ u/g(u) \mid u \in S \} \cup \{ 0\}$ is a hyperoval if and only if it is a Vandermonde set. 


\end{remark}

%
\begin{example} Vandermonde sets of the form $\{ u/g(u) \mid u \in S \} \cup \{ 0\}$ are not necessarily hyperovals either. For $q=16$, 
let 
$$	g(u) =  u^{16} + \omega u^{12} + \omega u^{11} + \omega u^6 + \omega u^5 + u + 1,
$$
where $\omega$ satisfies $\omega^3=1$. 	Let $H:= \{ u/g(u) \mid u \in S \} \cup \{ 0\}$. 
Then $\pi_i=0$, for $i=\{ 1, 3, 5, 7, 9, 11, 13\}$, but $\pi_{37} \ne 0$. 
This implies $H$ is a Vandermonde set but not a hyperoval. 
\end{example} 	
%
%
%
%
%

%
%
\begin{table}[h]
\centering

\begin{tabular}{|c|c|c|}
	\hline
	$q$  &   Elements in $\mathcal{D}$ & $\mid \mathcal{D} \mid$ \\  
	\hline
	4  &  1  & 1 \\
	\hline
	8 &	 1,3,5 & 3\\
	\hline
	16  & 1,3,5,7,9,11,13,37 & 8 \\
	\hline
	32   & 1, 3, 5, 7, 9, 11, 13, 15, 17, 19, 21, 23, 25, 27, 29, 69, 73, 77, 85, 89, 147 & 21  \\ 
	\hline 
	64	& 1, 3, 5, 7, 9, 11, 13, 15, 17, 19, 21, 23, 25, 27, 29, 31, 33, 35, 37, 39, 41, 43, & 55 \\
	&45, 47, 49, 51, 53, 55, 57, 59, 61, 133, 137, 141, 145, 149, 153,157, 165, 169, &\\
	&  173, 177, 181, 185, 275, 281, 283, 291, 297, 299, 307, 313, 409, 425, 661 &\\
	\hline
	128 &  1, 3, 5, 7, 9, 11, 13, 15, 17, 19, 21, 23, 25, 27, 29, 31, 33, 35, 37, 39, 41, 43, & 147 \\ 
	& 45, 47, 49, 51, 53, 55, 57, 59, 61, 63, 65, 67, 69, 71, 73, 75, 77, 79, 81, 83, 85, & \\
	& 87, 89, 91, 93, 95, 97, 99, 101, 103, 105, 107, 109, 111, 113, 115, 117, 119, 121, & \\
	& 123, 125, 261, 265, 269, 273, 277, 281, 285, 289, 293, 297, 301, 305, 309, 313,  &\\
	&  317, 325, 329, 333, 337, 341, 345, 349, 353, 357, 361, 365, 369, 373, 377, 529,&\\
	& 531, 537, 539, 547, 553, 555, 561, 563, 569, 571, 579, 585, 587, 593, 595, 601,&\\
	& 603, 611, 617, 619, 625, 627, 633, 785, 793, 809, 817, 825, 841, 849, 	857, 873, & \\
	& 881, 1093, 1095, 1107, 1109, 1111, 1123, 1125, 1127, 1139, 1141, 1301, 1317, & \\
	&  1333, 1365, 1381, 1587, 1619, 2341, 2349, 2381, 2405 &
	\\
	\hline
\end{tabular}

\caption{Set $\mathcal{D}$ in small fields}
\label{tabledset}
\end{table}

\section{The coefficients of $g$-functions and $\rho$-polynomials}
In this section we study the coefficients of $g$-functions for hyperovals. We start with the following. 
\begin{lemma} \label{lagrangemod} For $f, g \in K[x]$, $f(u) = g (u)$ for all $u \in S$ if and only
if $f (x) \equiv g (x) \pmod{x^{q+1} - 1}$.
\begin{proof} 
	By the division algorithm, there exist $r,h \in K[x]$ such that
	$$
	f (x) - g (x) = h (x) (x^{q+1}-1) +	r (x),
	$$
	where $\deg r(x) <q+1$. Then $f (u) = g (u)$ for all $u \in S$ if and only if
	$r (u) = 0$. This occurs if and only if $ r(x) = 0$,  if and only if  $f (x) \equiv g (x) \pmod{x^{q+1} - 1}$.
\end{proof}

\end{lemma}

For $1 \le k \le q$, let
$$
e_k(X_1, \dots, X_n) = \sum_{1 \le j_1 <\cdots<j_k \le n} X_{j_1} \cdots X_{j_k}
$$
be the elementary symmetric polynomial of degree $k$ for variables $X_1, \dots , X_n$.


\begin{lemma} \label{sympoly} $e_k(u_1, \dots, u_{q+1})=0$ for $1 \le k \le q$.
\begin{proof} To simplify the notation, we let $e_k=e_k(u_1, \dots, u_{q+1})$.
	We consider the following forms of the polynomial 
	$$
	f(x)=x^{q+1}-1 = \prod_{u \in S} (x-u) = x^{q+1}+e_1 x^{q}+e_2 x^{q-1} +\cdots +e_qx+e_{q+1}. 
	$$
	Then we have $e_k=0$ for $1 \le k \le q$.
\end{proof}
\end{lemma}

We now recall the Lagrange's Interpolation Formula (cp. \cite{lidl1997} Theorem 1.71). For $n \ge 0$, let $a_0,\dots ,a_n$ be $n + 1$ distinct elements of a field $F$, and let $b_0,\dots ,b_n$ be $n+1$ arbitrary elements of $F$. Then there exists exactly one polynomial $f \in F[x]$ of degree	$\le n$ such that $f(a_i)=b_i$ for $i = 0,\dots, n$. This polynomial is given by
$$
f(x)= \sum_{i=0}^{n} b_i \prod_{\substack{k=0 \\k \ne i}}^n \dfrac{x-a_k}{a_i-a_k}. 
$$

\begin{theorem} \label{thmcoefg} Let $H:= \{ u /g(u) \mid u \in S \} \cup \{ 0\}$. Then $H$ is a hyperoval if and only if one of the following holds.
\begin{enumerate}
	\item The coefficient of $x^{2i-k+q+1}$ in $g^{q-1-k}(x) \pmod{x^{q+1}-1}$ is zero,  for each pair $(i,k) \in \mathscr{M}$.
	\item  The coefficient of $x^{k-2i}$ in $g^{q-1-k}(x) \pmod{x^{q+1}-1}$ is zero,  for each pair $(i,k) \in \mathscr{M}$.
\end{enumerate}	
\begin{proof} Since the coefficients of $x^{k-2i}$ and $x^{2i-k+q+1}$ of $g^{q-1-k}(x)$ are conjugate, it is sufficient to prove that $H$ is a hyperoval if and only if the condition in part 1) holds. 
	We can  assume $g \in K[x]$ such that $g(u) \in F$ for $u \in S$. From the Lagrange's Interpolation Formula and  Lemma \ref{lagrangemod}, for each $k$, modulo  $(x^{q+1}-1)$, $g^{q-1-k}(x)$ has the form
	$$
	g^{q-1-k}(x) \equiv \sum_{u \in S} g^{q-1-k}(u) \prod_{\substack{v \in S \\v \ne u}} \dfrac{x-v}{u-v}. 
	$$
	Using Lemma \ref{sympoly}, for $u \in S$, we have
	\begin{align*}
	\prod_{\substack{v\in S\\ v \ne u}}(x-v) =
	x^q+ux^{q-1}+u^2x^{q-2}+ \cdots + u^q.
	\end{align*}
	Also, 
	\begin{align*}
	\prod_{\substack{v\in S\\ v \ne u}}(u-v)= u^q.
	\end{align*}
	
	Then,
	\begin{align*}
	g^{q-1-k}(x)&\equiv\sum_{u \in S} \left[ ug^{q-1-k}(u) (x^q+ux^{q-1}+ \cdots + u^q) \right]\\
	&= \sum_{i=1}^{q+1}\left(\sum_{u \in S} u^ig^{q-1-k}(u) \right) x^{q+1-i}	\\
	&=\sum_{i=0}^{q}\left(\sum_{u \in S} u^{-i}g^{q-1-k}(u) \right) x^{i}.
	\end{align*}
	By Theorem \ref{thmhovalD2}, $H$ is hyperoval if and only if 
	$$
	\sum_{u \in S} \left(\dfrac{u}{g(u)}\right)^d
	= 
	\sum_{u \in S} \left(\dfrac{u}{g(u)}\right)^{iq+k-i}
	= 
	\sum_{u \in S} u^{k-2i}g^{q-1-k}(u)
	=0,
	$$
	for each pair $(i,k) \in \mathscr{M}$. 
	This occurs if and only if the coefficient of $x^{2i-k+q+1}$ in $g^{q-1-k}(x) \mod (x^{q+1}-1)$ is zero. \qedhere

\end{proof}
\end{theorem}
Using the relationship $\rho(u)=1/g(u)$ we  obtain the following. 
\begin{theorem} \label{thmcoeff}  Let $H:= \{ u \rho(u) \mid u \in S \} \cup \{ 0\}$. Then $H$ is a hyperoval if and only if one of the following holds.
\begin{enumerate}
	\item The  coefficient of $x^{2i-k+q+1}$ in $\rho(x)^k \pmod{x^{q+1}-1}$ is zero, for each pair $(i,k) \in \mathscr{M}$.
	\item  The   coefficient of $x^{k-2i}$ in $\rho(x)^k  \pmod{x^{q+1}-1}$ is zero, for each pair $(i,k) \in \mathscr{M}$.
\end{enumerate}	
\end{theorem}

\begin{corollary} If $H:= \{ u/ g(u) \mid u \in S \} \cup \{ 0\}$ is a hyperoval, then the coefficient of $x^t$ in $g(x)$ is zero,
for $t\equiv 2 \pmod 4$ or $t\equiv  3 \pmod 4$.

\begin{proof} Assume $H:= \{ u/ g(u) \mid u \in S \} \cup \{ 0\}$ is a hyperoval. We consider the special case $k=q-2$. We have
	$$
	\mathscr{I}:=\{ i \mid (i,k) \in \mathscr{M}  \}=\{ i \mid  0\le i \le \lfloor(k-1)/2\rfloor, 
	i \text{ even}\}.
	$$
	This implies 
	$$
	\{ k-2i \mid  i \in \mathscr{I}\} = \{ t \mid 1<t<q,t \equiv 2 \pmod 4\}.
	$$
	
	By Theorem \ref{thmcoefg}, the coefficient of $x^{t}$ is zero in $g(x)$ for  $ t \equiv 2 \pmod 4$.  The proof now follows from the fact that the coefficients of $x^t$ and $x^{q+1-t}$ are conjugate. 
\end{proof}
\end{corollary}

\section{Gram matrices}

Following \cite{fisher2006}, consider an element $\mathbf{i} \in K$ with property $T(\mathbf{i}) = \mathbf{i}+\mathbf{i}^q = 1$. Then $K = F(\mathbf{i})$ and $\mathbf{i}$ is a root of a quadratic equation
$$
z^2 + z + \delta = 0,
$$
where $\delta = N(\mathbf{i}) \in F$. Any element $z \in K$ can be represented as $z = x + y\mathbf{i}$, where $x, y \in F$. For $z = x + y\mathbf{i}$ we have $x = \langle \mathbf{i}, z\rangle$, and $y = \langle 1, z \rangle$.

Note that if $m$ is odd then one can choose $\mathbf{i} = \omega, \omega^2 + \omega + 1 = 0$. So if $w \in K$ is a
generator of $S$ then we can take $\mathbf{i} = \omega = w^{(q+1)/3}.$
From Theorem  \ref{thmbracket} we have the following. 

\begin{theorem} \label{hovalmatrix} Let $H:=\{ y_i \mid 1\le i \le q+1\} \cup \{ 0\}$. Assume no two points $y_i, y_j$ are on the same line $\langle v, x\rangle = 0$ for all $v \in S$. Then $H$ is a hyperoval if and only if 
$$
\sum_{j=1}^{q+1} \langle y_i, y_j \rangle^k = 0,
$$
for each $y_i$ and $1 \le k \le q$. 
\end{theorem}
For each $1\le k\le q$, let $M_H(k)$ be the $(q+1) \times (q+1)$ matrix whose entries are $(M_H(k))_{i,j}:= \langle y_i, y_j \rangle^k$.
By Theorem \ref{hovalmatrix}, sums of elements in any column and any row of $M_H(k)$ are equal to $0$. 

In the following we consider additional properties of the matrix $M_H(1)$, which we will denote by $M_H$. 

\begin{proposition} Let $H:=\{ y_i \mid 1\le i \le q+1\} \cup \{ 0\}$ be a hyperoval.  For each $i$, let $c_i,s_i$ be elements in $F$ such that  $y_i=c_i+s_i\mathbf{i}$.  Without loss of generality, assume $y_q=1$ and $y_{q+1}= s_{q+1} \mathbf{i}$.  

Let $M_H$ be the $(q+1) \times (q+1)$ matrix whose entries are $(M_H)_{i,j}:= \langle y_i, y_j \rangle$. The matrix $M_H$ has the following properties.
\begin{enumerate}
	
	\item  $M_H$ is symmetric and its diagonal entries are zeros. The trace of  $M_H$ is zero. 
	
	\item $M_H$ has an eigenvalue $\mu =0$ of multiplicity $q-1$. The corresponding eigenvectors are 
	$$
	v_i=[0, \dots , 0, \stackrel{i}{1},0,\dots,c_i,s_i/s_{q+1}]^T.
	$$
	where $1 \le i \le q-1$.
	
	Also, the determinant of $M_H$ is zero.
	
	\item  The characteristic polynomial of $M_H$ is $P(x)=x^{q-1}(x+\mu_0)^2$. In particular, $M_H$ has exactly one non-zero eigenvalue $\mu=\mu_0 \in F$ of  multiplicity $2$. 
	
	\item The minimal polynomial of $M_H$ is  $Q(x)=x(x+\mu_0)$.
	
\end{enumerate}

\begin{proof}
	
	\begin{enumerate}
		
		\item These follow from properties of the bilinear form $\langle \cdot,\cdot \rangle$.
		
		\item From Theorem \ref{hovalmatrix}, it can be calculated that $M_Hv_i=0$. Since the set of vectors $\{ v_i\}$ is linearly independent, $\mu=0$ is an eigenvalue of $M_H$ with multiplicity $q-1$ and $v_i$'s are the corresponding eigenvectors. It follows that the determinant of $M_H$ is zero. 
		
		\item 	Since $\mu =0$ is an eigenvalue of multiplicity $q-1$, the characteristic polynomial of $M_H$ has the form
		$$
		P(x)=x^{q-1}(x^2+ax+\mu_0^2).
		$$
		Also since $tr(M_H)=0$, the coefficient $a$ of $x^q$ is zero. Hence
		$
		P(x)=x^{q-1}(x+\mu_0)^2.
		$
		This implies $\mu=\mu_0$ is the only non-zero eigenvalue of $M_H$ with multiplicity $2$. 
		\item 
		The minimal  polynomial $Q(x)$ of $M_H$ divides $P(x)$, so it has the form
		$$
		Q(x)=x^s(x+\mu_0)^t.
		$$  
		Since $M_H$ is diagonalizable, $s=t=1$ and so  
		$
		Q(x)=x(x+\mu_0). \qedhere
		$
	\end{enumerate}	
\end{proof}
\end{proposition} 

{\bf Acknowledgments}

\medskip

This work was supported by grant 31S366.

 \printbibliography

\end{document}